\documentclass[12pt]{article}
\advance \topmargin by -\headheight
\advance \topmargin by -\headsep
     
\evensidemargin \oddsidemargin
\newcommand\be{\begin{equation}}
\newcommand\ee{\end{equation}}
\newcommand\bea{\begin{eqnarray}}
\newcommand\eea{\end{eqnarray}}

\begin{document}
\title{\bf Stirling's formula derived simply}
\author{Joseph B. Keller\\Departments of Mathematics and Mechanical Engineering\\Stanford University, Stanford, CA  94305-2125\\email: keller@math.stanford.edu\\
and\\
Jean-Marc Vanden-Broeck\\
Department of Mathematics\\University College London\\ London WC1E 6BT, United Kingdom\\email:j.vanden-broeck@ucl.ac.uk}

\date{November 26, 2007}
\maketitle

\begin{abstract}
Stirling's formula, the asymptotic expansion of $n!$ for $n$ large, or of $\Gamma(z)$ for $z\to \infty$, is derived directly from the recursion equation$\Gamma(z+1) =z \Gamma(s)$ and the normalization condition $\Gamma \left(\frac{1}{2}\right) =\sqrt{\pi}$.
\end{abstract}

\newpage

Stirling's formula is an asymptotic expansion of the gamma function $\Gamma(z)$, valid as $z\to \infty$ with  $| \arg z|<\pi$.  Since $\Gamma (z)$ is the analytic extension of $(n-1)!$, it satisfies the recursion equation 
\be
\Gamma (z+1) =z\Gamma (z).
\label{1}
\ee
We shall derive Stirling's  formula from (\ref{1}) and the condition
\be
\Gamma \left(\frac{1}{2}\right)=\sqrt \pi,
\label{2}
\ee
which makes the solution unique.

Taking the logarithm of each side of (\ref{1}), and writing $f(z) =\log \Gamma (z)$, yields
\be
f(z+1)- f(z) =\log z.
\label{3}
\ee
Replacing $f(z+1)$ in (\ref{3}) by its  Taylor series around $z$  leads to 
\be
\sum\limits^\infty_{n=1} \: \frac{f^{(n)}(z)}{n!} =\log z.
\label{4}
\ee
We integrate (\ref{4}), introducing a constant of integration $\log C$, and we call the right side $f_0(z)$:
\be
\sum\limits^\infty_{n=1} \:\frac{f^{(n-1)}(z)}{n!} =z\log z-z+\log C \equiv f_0 (z).
\label{5}
\ee
Keeping the undifferentiated term with $n=1$ on the left, and moving the other terms to the right, gives 
\be
f(z) =f_0 (z) -\sum\limits^\infty_{n=2}  \: \frac{f^{(n-1)} (z)}{n!}.
\label{6}
\ee

Equation (\ref{6}) is a differential equation for $f$ of infinite order, which we   solve by successive substitution. First we approximate $f$ by $f_0$.  Then we use $f_0$ in place of $f$ in the sum to get a second approximation to $f$.  Then we substitute that for $f$ in the sum, and so on.  This will lead to an expansion for $f$ of the form
\be
f(z) \sim\sum\limits^\infty_{j=0} c_j \, f^{(j)}_0 (z).
\label{7}
\ee
We expect (\ref{7}) to be an asymptotic expansion of $f(z)$, valid as $z\to \infty$.  Therefore we have used the sign $\sim$ of asymptotic equality.

To determine the coefficients $c_j$ in (\ref{7}) we substitute (\ref{7}) into (\ref{6}) and equate coefficients of $f^{(j)}_0$.  This yields
\be
c_0 =1, \qquad c_j =-\sum\limits^{j-1}_{k=0} \: \frac{c_k}{(j-k+1)!} , \qquad j\geq 1.
\label{8}
\ee
Equation (\ref{8}) determines the $c_j$ successively, starting with $c_0 =1$,  $c_1 = -1/2$, etc. To solve it in general we set $c_j =B_j /j!$ in (\ref{8}).  The resulting recursion equation for $B_j$ is exactly that satisfied by the Bernoulli numbers [eq.\ (6.79), p.\ 270 of Graham, Knuth and Patashnik \cite{GKP}].   Thus  $B_j$ is the $j$-th Bernoulli number. 

The   expansion (\ref{7}) of $f$ involves the derivatives of $f_0$.    Since $f^{(1)}_0 =\log z $ it follows that $f^{(j)}_0 =(-1)^j (j-2) !/ z^{j-1}$ for $j\geq 2$.   For $j\geq 2$ with $j$ odd, 
$B_j =0$.   Therefore by setting $j=2n$ for $j\geq 2$, and using the definition $f(z)=\log\Gamma (z)$,  we can  write (\ref{7}) as follows:
\be
\log \Gamma (z) \sim \left( z-\frac{1}{2}\right) \log z -z +\log C+\sum\limits^\infty_{n=1} \:\frac{ (2n-2)! B_{2 n}}{(2n)! \, z^{2n-1}}.
\label{9}
\ee
Exponentiating both sides of (\ref{9}) gives
\be
\Gamma (z) \sim C z^{z-\frac{1}{2}} e^{-z} \exp \left[ \sum\limits^\infty_{n=1} \: \frac{(2n-2)! \; B_{2n}}{(2n)! \; z^{2n-1}}\right].
\label{10}
\ee

To find $C$ we first prove that
\be
\Gamma \left(n+\frac{1}{2}\right) =\frac{\Gamma (2n+1) \sqrt \pi}{2^{2n} \,\Gamma (n+1)}.
\label{11}
\ee 
For $n=0$ this reduces to (\ref{2}) , so it is true.  It follows for $z=n$ by induction on $n$ using (\ref{1}).  When we use  the leading term in (\ref{10}) for $\Gamma(n)$ in (\ref{11}),  we   find that 
\be 
C =\sqrt{2\pi}.
\label{12}
\ee
This completes the determination of the expansions (\ref{9}) and (\ref{10}) of $\log \Gamma$ and of $\Gamma$.  With $C=\sqrt{2\pi}$, (\ref{9}) agrees with eq.\ (6.1.40) of Abramowitz and Stegun \cite{AS}.  When the exponential in (\ref{10}) is   expanded in a power series, it agrees with eq.\ (6.1.37) of Abramowitz and Stegun \cite{AS}, which is Stirling's formula.

\end{document}